\documentclass[12pt]{amsart}
\usepackage{amsmath,amssymb,amscd}
\usepackage[mathscr]{eucal}
\newcommand{\Spec}{{\mathrm{Spec}\, }}

\newcommand{\bZ}{{\mathbb Z}}
\newcommand{\bR}{{\mathbb R}}
\newcommand{\bC}{{\mathbb C}}
\newcommand{\bN}{{\mathbb N}}
\newcommand{\cA}{{\mathcal A}}
\newcommand{\cB}{{\mathscr B}}
\newcommand{\cC}{{\mathscr C}}

\newcommand{\cF}{{\mathcal F}}
\newcommand{\cG}{{\mathcal G}}

\newcommand{\cI}{{\mathcal I}}
\newcommand{\cJ}{{\mathcal J}}

\newcommand{\cL}{{\mathcal L}}

\newcommand{\cM}{{\mathcal M}}
\newcommand{\cN}{{\mathcal N}}
\newcommand{\cO}{{\mathcal O}}
\newcommand{\cP}{{\mathscr P}}

\renewcommand{\sp}[1]{{\mathrm{Spec}}}

\theoremstyle{plain}
\newtheorem{thm}{Theorem}[section]

\newtheorem{pro}[thm]{Proposition}

\theoremstyle{definition}

\newtheorem{voi}[thm]{}

\begin{document}
\bigskip 
\title{Notes on Nash modification}
\author{Augusto Nobile}


\address{Louisiana State University \\
Department of Mathematics \\
Baton Rouge, LA 70803, USA}



%
\email{nobilba@gmail.com}
\maketitle
\begin{abstract}

The Nash blowing-up (or modification) of an algebraic variety $X$ is a canonical process  that produces a proper, birational morphism $\pi : X' \to X$  of varieties. It is expected that the singularities of $X'$ will be better than those of $X$.  In the mid-1970's, it was proved that in characteristic zero, ${\pi}$ is an isomorphism if and only if $X$ is nonsingular, which is false in positive characteristic. The focus of this article is on several subsequent studies on this subject.
 Topics covered  include: (a) the extension of the mentioned theorem to the case where $X$ is normal, in any characteristic, (b) the introduction and study of  Nash modifications of higher order, (c) the case where the variety $X$ is toric, where more precise results can be obtained and (d) desingularization properties of the Nash process.
 
\end{abstract}

\section*{Introduction}
\label{S:intro} 

The so-called Nash blowing-up is a procedure to modify an algebraic variety $X$ into another $\cN(X)$, equipped with a proper birational morphism $\pi : \cN(X) \to X$, which is an isomorphism over the set
$X_{reg}$ of regular, or non singular points of $X$. It is expected that somehow the singular locus of $X$ will improve in the process.

Intuitively, the idea is to replace each singular point by the limits of tangent spaces at nearby regular points. 

The term ``blowing up'' is perhaps inaccurate, since in general this is not the blowing up of a specific $\cO _X$-ideal. Rather, it is a canonical process that does not involve the choice of a center on $X$. Perhaps it should be called ``Nash modification or transformation'', instead. Generally we shall use the expression ``Nash modification'', to my knowledge a term proposed by Teissier.

It seems that the first place where the expression ``Nash blowing-up'' appears in print is the article \cite{N}, written by the author of this note in the mid 1970's. I heard the term (and its description) in a talk by John Mather, in his seminar at Harvard University around 1971. The seminar dealt primarily with Thom's theory of stratified spaces and related matters. It seems that the procedure was in  the ``folklore''. 
Much later, it was noticed that the operation appears in \cite{Se}, a paper by Semple, published in 1954. It is not clear whether Semple created the method, or it was already around in the right circles. Anyway, Semple's contribution was ignored until 2014, when 
 Gonz{\'a}lez Perez and Teissier mentioned in \cite{GT}. In this article they suggest, reasonably, to call it  the Semple-Nash 
modification. 
 But in the present article we shall use the expression  {\it Nash modification}. 
 
 It seems that Nash did not publish anything where the procedure is used. 

In \cite{N} it is proved that working over $\bC$, the morphism $ \pi: \cN(X) \to X$ is an isomorphism if and only if $X$ is nonsingular. From this it easily follows that the iteration of the process desingularizes a curve. On the contrary, for the plane curve $X:y^2-x^3=0$, working over a field of characteristic two, the Nash modification is an isomorphism, although of course $X$ is singular. So, a natural question would be whether  
working in characteristic zero, the (finite) iteration of the Nash modification desingularizes a variety of dimension $>1$.

Not surprisingly, the best results are in dimension 2 (surfaces). In his thesis \cite{Re} written in the late 1970's (under the direction of Morrow, seemingly unpublished), Rebassoo shows that indeed the iteration of the Nash modification resolves the singularities of surfaces of a special type. They are defined in ${\mathbb A}^3 $ by certain equations and (this is very special) the successive transforms are again locally hypersurfaces in ${\mathbb A}^3$, whose equations can be controlled.

Later other results, more general because they apply to any surface in characteristic zero, were obtained. However, they do not deal strictly with the Nash modification, but with the normalized one, i.e, where it is followed by normalization. Using this process, Hironaka made important contributions to the problem (\cite{Hi}) and, based on this work, Spivakovsky proved that  iteration of normalized Nash blowing-ups desingularizes surfaces in characteristic zero (\cite{S}). 

In positive characteristic it seems that what is essential in the mentioned example in \cite{N} is the fact that the considered curve, being singular, is not normal. Indeed, in the recent article \cite{DNB} the authors show that if $X$ is a normal variety of any dimension, the Nash morphism $\pi : \cN(X) \to X$ is an isomorphism if and only if $X$ is non singular, even if the base field has positive characteristic.  So, Spivakovsky's result might be true in positive characteristic, as well as its analog for dimension $> 2$.

There are generalizations of the construction of the Nash transformation. For instance, in \cite{OZ} the authors study a modification of a variety $X$ relative to a given coherent $O_{X}$-module $\cF$.  In particular, taking as  $\cF$   a suitable sheaf, we get the so-called Nash modification of order $n$, $\pi_n:\cN_n(X) \to X$. For $n=1$ we find again the morphism $\pi$ mentioned above. Of course, again we may consider the ``normalized'' version  $\pi^{-}_n: \cN_n^{-}(X) \to X$, that is where 
$\cN(X)$ is substituted by its normalization $\cN_n^{-}(X) $.  

In \cite{Ya}, Yasuda carefully studies the modification $\pi_n$, in particular he gives some equivalent versions of the definition.  He also proves that for a curve $X$ (in characteristic zero) $\pi_n :\cN_n(X) \to X$ is a desingularizarion of $X$, for $n$ large enough. 
 If, as he speculates, a similar statement were valid in arbitrary dimension, we would have a remarkable result: a one step canonical desingularization of any variety. Unfortunately, this is not true. In \cite{TY} Toh-Yama gave an example of a surface $X$ for which 
$\cN_n^{-}(X)$ is singular for any natural $n$. 

If we focus our attention to toric varieties, more general results may be obtained about the Nash morphism $\pi$ and also for $\pi_n$, even over fields of positive characteristic, see \cite{D} and \cite{DNB}. The example of Toh-Yama in \cite{TY} uses techniques from Toric Geometry.

In this article we intend to discuss in greater detail the topics mentioned in this introduction, with the exception of the work of Hironaka and Spivakovsky (\cite{Hi}, \cite{S}). In general we won't give detailed proofs, but we'll try to indicate the main ideas and techniques behind them.

  \section{Defining the Nash modification}
\label{S:defi}

\begin{voi}
\label{d1}
For us, ``variety'' means an integral scheme of finite type over an algebraically closed field $k$. As is well known, in that case closed points are dense, and they play an important role (see \cite{Ha},  ). In many instances we may restrict ourselves to them. Hence, in general, for us ``point'' will mean a closed point. A nonsingular point will also be called a regular point, $X_{reg}$ denotes the (open) set of regular points of $X$.
\end{voi}

\begin{voi}
\label{d2}
 The most elementary way to  introduce the Nash modification of a $d$-dimensional variety $X$ uses the classical Grassmanian $G(d,N)$: the variety parametrizing vector subspaces of dimension $d$ of $k^N$. Let $X$ be a $d$-dimensional closed subvariety of $\mathbb A ^N$. If $x \in X$ is a regular point, it has a tangent space $T(X,x)$ of dimension $d$, which can be identified to a linear subspace of $k^N$, i.e to a point of $G(d,N)$. Thus we have a morphism $\phi:X_{reg} \to X \times G(d,N)$, where 
$\phi (x) =(x, T(X,x))$. Let $X'$ be the closure of the image of $\phi$ in $G(d,N)$. The first projection induces a proper, birational morphism $\pi: X' \to X$ which is an isomorphism above $X_{reg}$. Then we define $X':= \cN(X)$, the Nash modification. 

It is easy to show that the morphism $\pi$ is independent of the embedding of $X$ into an affine space, as a consequence this operation extends to an abstract variety. With this approach it is clear that in this process a singular point of $X$ is replaced by the limiting positions of tangent spaces taken at nearby regular points.

\end{voi}

\begin{voi}
\label{d3}
 A slightly different way to phrase the construction above is as follows. Let $X \subset {\mathbb A}^N$ be our variety. If $x \in X$ is regular and $(\cO _{X,x},\cM _x)$ its local ring, recall that we may identify 
$T(X,x)$ and the dual ${(\cM_x/\cM_x ^2)}^{\star}$ of the $k$-vector space ${(\cM_x/\cM_x ^2)}$.
This ``globalizes'': if 
$R=\Gamma(X, \cO _X)$ is the affine algebra of $X$, $I$ is the kernel of the multiplication homomorphism 
$R \otimes _k R \to R$ and we consider $I/I^2$ (an $R$-module, via
 $r(a \otimes b)= ra \otimes b$, usually called  the {\it conormal module}), then
 ${(I/I^2)}_x \otimes _R k = \cM_x/{\cM_x}^2$, whence
 ${({(I/I^2)}_x \otimes_R k)}^{\star}=T(X,x) \subset T(\mathbb A ^N,x)$, where $V^{\star}$ denotes the dual of a vector space $V$.
 
 As before, with $U=X_{reg}$, we have a map
 $\phi : U \to G(d,N)$ , $u \to (u, T(X,u))$. The closure of the image of $\phi$ is $\cN(X)$.
 \end{voi}
 
 \begin{voi}
 \label{d4}
  This approach suggests a generalization: instead of the $R$-module $I/I^2$ we could take $I/I^{n+1}$, for any integer $n \ge 2$. Let 
 $$T^n(X,x):=  {(( I_x/I_x ^{n+1}) \otimes_R k)} ^{\star} ~.$$
 If $x$ is regular, the dimension of this vector space is $\delta=\binom{d+n}{d}-1$. We have an inclusion  of 
 $T^n(X,x)$ into  $ T^n(\mathbb A ^N,x)$,whose dimension  is $\delta'=\binom{N+n}{N}-1$. As before we may consider the corresponding map $\phi_n:U \to U \times G(\delta,\delta')$.  The closure of $Im(\phi_n)$, denoted by  $\cN _n(X)$, is the $n$-th Nash modification of $X$. We have a natural projection $\pi _n : \cN_n(X) \to X$, which is proper and an isomorphism over $U=X_{reg}$.
 \end{voi}

\section{The basic results}
\label{S:resu}

Probably the most basic result related to the Nash modification is:

\begin{thm}
\label{basic}
 Let $X$ be an algebraic variety over a field $k$, algebraically closed of zero characteristic , 
$\pi: \cN(X) \to X$ the Nash modification.  Then $X$ is regular if and only if $\pi$ is an isomorphism.
\end{thm}

\bigskip

Seemingly, the first proof of this result appears in \cite{N}. There one works over $\mathbb C$, the complex numbers, and the variety $X$ is regarded as a complex analytic variety, mainly methods from local analytic geometry are used.  In the proof, one reduces the general case to that where $X$ is one dimensional. The modification is understood in the sense of \ref{d2}.  We won't review that proof here.

\bigskip

Another, more algebraic (and simpler) proof is due to Teissier (\cite{T}, pages 585-587). In it, a different (equivalent) approach to the definition of the Nash process is used. In fact, if $X$ be a variety of pure dimension $d$, over an algebraically closed field $k$ and $\Omega = \Omega _{X/k}$ its sheaf of differentials, then there is birational, projective morphism $\pi : X_1 \to X$ and an exact sequence of $\cO _{X_1}   $-modules 
$$(1)  \quad {\pi}^{\star} (\Omega) \to \cL \to 0 \,$$
 with $\cL$ locally free of rank $d$, universal in the sense that if $h:Z \to X$ is another morphism such there is an exact sequence of $\cO _{Z}   $-modules 
${h}^{\star} (\Omega) \to \cG \to 0 $ 
with $\cG$ locally free of rank $d$, then there is a unique morphism $q:Z \to X_1$ such that ${\pi}q=h$ and $q^{\star}(\cL)=\cG$. This pair $(X_1,\pi)$ is unique up to canonical isomorphism. This assertion is an easy consequence of 
Grothendieck's theory of the Grassmanian, which will be recalled in Section  \ref{S:OZ}. 

It turns out that $X_1$ is isomorphic to  the Nash modification $\cN(X) $ of $X$ introduced in \ref{d2}. This is a consequence of the following remark.   
 If $x \in X$ is a closed point, the fiber $\pi^{-1}(x)$ parametrizes the $d$-dimensional quotients of the vector space $\Omega_x$. Dualizing, each of these quotients becomes a $d$-dimensional subspace of $T_{X,x}$. If $x$ is a regular point, we get  precisely the Zariski tangent space to $X$ at $x$. 
 
Now assume $\pi:X_1:=\cN(X) \to X$ is an isomorphism. Then we may identify $\Omega$ and $\pi^{\star}\Omega$. Since $\pi^{\star}\Omega$
has a locally free quotient (namely $\cL$ in (1)), $\Omega$ also has one. But we have the following result.

\begin{pro} 
\label{nons}
Let $X$ be an $d$-dimensional algebraic variety over an algebraically closed field of characteristic zero. Assume $\Omega_X $has a locally free quotient of rank $d$. Then $X$ is nonsingular.
\end{pro}

Now assume $\pi:X_1:=\cN(X) \to X$ is an isomorphism. Then we may identify $\Omega$ and $\pi^{\star}\Omega$. Since $\pi^{\star}\Omega$
has a locally free quotient of constant rank $d$ (namely $\cL$), $\Omega$ also has one. But we have the following result.

\begin{pro} 
\label{nons}
Let $X$ be an $d$-dimensional algebraic variety over an algebraically closed field of characteristic zero. Assume $\Omega_X $has a locally free quotient of rank $d$. Then $X$ is nonsingular.
\end{pro}

We sketch a proof of this result. We are easily reduced to the following more algebraic statement. Let $x \in X$ be a closed point, $A$ the completion of $\cO_{X,x}$. Assume we have an exact sequence
$$ (2) \qquad \Omega_{A/k} \to A^d \to 0$$ 
(where the module of differentials is in the formal, or analytic, sense). Then $A$ is regular. 

This is proved by induction on $d$ (the case $d=0$ is trivial, because then $A=k$). For the inductive step, one uses the following result, whose proof,
 using the formal exponential, is due to Zariski.

Let $A$ be as above, with maximal ideal $M$. Assume $D:A \to A$ is a $k$-derivation such that $D(A) \nsubseteq M$. Then, there is an element $x \in A$ and a local subring $A_1 \subset A$ such that $A=A_1[[x]]$, where $\dim (A_1) = d-1$.

To conclude the inductive step, given the sequence (2), we may choose a homomorphism $s:A^d \to \Omega_{A/k}$ such that $s(e_1)=1$ and $s(e_j)=0$ if $j>1$, where $e_1, \ldots, e_d$ is the standard basis of $A^d$. From this,  we may  get  $x \in A$ and a derivation $D$ such that $D(x)=1$. Using 
 Zariski's result we may write $A=A_1[[x]]$. But then we get  a sequence similar to (2), but with $A$ substituted by $A_1$ and $d$ by $d-1=\dim (A_1)$. Now we may use induction.

\section{The results of Anna Oneto and Elsa Zatini}
\label{S:OZ}

The results of \ref{d3} may be presented in a more intrinsic and general way by using Groethendieck's version of the Grassmanian, already mentioned in Section \ref{S:resu}. This is done in the article \cite{OZ}, we explain some of its contents.

Recall that if $S$ is a noetherian scheme, $\cF$ a quasi-coherent $\cO_S$-module and $n \in \bN$,  we may define a contravariant functor 
$$Grass_{\cF,n} : \{ S {\mathrm{-schemes}} \} \to \{\mathrm{Sets}  \}$$ 
by: given an $S$-scheme $f:T \to S$, $Grass_{\cF,n}(T)$ is the set of $n$-quotients of $f^{\star}{\cF}$.   Here, by an $n$-quotient of an $\cO_T$-module $\cG$, we mean a surjective homomorphism of $\cO_T$-modules $q:\cG \to \mathcal Q$, with $\mathcal Q$ locally free, of constant rank $n$, we identify $q$ and  
$q':\cG \to \mathcal Q'$ if there is an isomorphism $\eta : \mathcal Q \to \mathcal Q' $ such that $q \eta = q'$.

This functor is represented by an $S$-scheme $g: \underline{Grass}_{\cF,n} \to S$, where the morphism $g$ is proper, together with a universal $n$-quotient   
$u:g^{\star}(\cF) \to {\mathcal  Q} _{\underline {    Grass}_{\cF,n}     } $ (\cite{GD}, Thm. 9.7.4).

 From now on we shall assume that $\cF$ is coherent and there is a dense open set $U \subset S$ such that the restriction $\cF_{|U}$ is locally free, of constant rank $d$, we shall write $G:={\underline {Grass}}_{\cF,d}$.  
 Then, since the identity of $\cF_{|U}$ is a $d$-quotient, by the fact that $(G,u)$ represents $Grass_{\cF,d}$, we get a unique morphism 
 $\sigma : U \to G$ of $S$-schemes such that 
 $\cF _{|U}={\sigma}^{\star}(\mathcal Q_{G})$. 
 
 Let $S^N(\cF)$ be the closure (in $G$) of $Im(\sigma)$, then there is morphism 
 $$v:S^N(\cF) \to S$$
 induced by $g:G \to S$.  The morphism $v$ is projective (because  $\cF$ is coherent).      
 The pair $(S^N(\cF), v)$ is the Nash modification of $S$ relative to $\cF$.  
 Alternatively we'll write $S^N(\cF)= \cN(S, \cF)$ in this case.

 The scheme $\cN(S, \cF)$ has the following universal property (inherited from that of $G$): 
 if $h:T \to S$ is such that there is a $d$-quotient $q: h^{\star}(\cF) \to \mathcal Q$, then there is a unique morphism $\rho:T \to \cN(S,\cF)$ such that $q$ is induced by 
 $u: g^{\star}(\cF) \to {\mathcal Q}_{G}$, satisfying the equality 
 $v \rho = h$.

 It can be proved that $S^N(\cF) \cong S^N(\bigwedge ^d \cF)$.
 
 An interesting special case is that where $S=X$, a $d$-dimensional irreducible algebraic variety over an algebraically closed field $k$. Let 
 $\Omega _X=\Omega _{X/k}$ be the sheaf of differentials on $X$. If $U=X_{reg}$, then  $\Omega_X$ restricted to $U$ is locally free, of rank $d$. Then,
 $ \cN(X,\Omega_X)  = X^N(\Omega_X)=\cN(X)$, the Nash modification in the sense of \ref{d2}. This is Teissier's approach in Section \ref{S:resu}. 
 
 The article \cite{OZ} contains other interesting results. For instance, $v:\cN(S,\cF) \to S$ is compared to $S' \to S$, the blowing up of $S$ with respect to the $d$-th Fitting ideal of $\cF$ ($d=\dim S$).  It is shown that there is a morphism of $S$-schemes 
 $u: S' \to \cN(S,\cF)$, a necessary and sufficient condition in order that $u$ be an isomorphism is that the projective dimension of 
 $v^{\star} \cF$ be $ \le 1$.
 
 They also prove that always the Nash modification $v:\cN(S,\cF) \to S$ is isomorphic to the blowing up of $S$ with respect to a fractionary ideal $\mathcal E \subset \mathcal R$, the sheaf of total rings of fractions of $\cO_S$. This ideal $\mathcal E$ may be explicitly described: up to some identifications, we get a homomorphism 
 $\phi : \wedge ^d \cF \to \mathcal R$, 
 and  $\mathcal E =  \phi (\wedge ^d \cF)$. Using this description of the Nash transformation, they find an alternative proof of Theorem \ref{basic}.
 
 Other results involving Jacobian ideals and geometrically linked varieties are also obtained in the article just reviewed.

\section{Yasuda's work}
\label{S:Yw}
In the paper \cite{Ya}, Yasuda presents a thorough discussion of the notion of Nash modification of orden $n$, as well as some fundamental results and  interesting conjectures. One works with algebraic varieties over an algebraically closed field $k$, as usual ``point'' generally means closed, or $k$-rational, point.

In \cite{Ya}, the basic proposed definition of  $\pi_n:\cN_n(X) \to X$, the $n$-th Nash modification of the $d$-dimensional variety $X$, uses a Hilbert scheme. Namely, if 
$({\cO}_{X,x},{\cM}_x)$ is the local ring of $X$ at a regular closed point  $x \in X$, then 
 ${\cO}_{X,x}/{\cM_x}^{n+1}$ is a $k$-vector space of dimension $\delta=\binom{d+n}{n}$.  Consider the $X$-scheme ${\mathrm H}_{\delta}(X)$,  parametrizing sets of $\delta$ points of $X$ and their specializations,    in particular subschemes $Z_x \subset X$ concentrated at a single regular point $x \in X$,  the only nontrivial stalk of $\cO_{Z_x}$ being  
${\cO}_{X,x}/{\cM_x}^{n+1}$.

 We have a morphism $\sigma_n:X_{reg} \to {\mathrm H}_{\delta}(X)$, $x \to [Z_x]$, with $[Z_x]$ the point of ${\mathrm H}_{\delta}(X)$ corresponding to $Z_x$.  The closure of the graph of $\sigma _n$ in $X \times {\mathrm H}_{\delta}(X)$ is, by definition, $\cN_n(X)$. Since ${\mathrm H}_{\delta}(X)
$ is an $X$-scheme, there is an induced projection $\pi_n:\cN_n(X) \to X$, which is a proper, birational morphism (an isomorphism over $X_{reg}$). This modification coincides with that described in \ref{d4}. 

Yasuda proves that in characteristic zero, if $X$ is a curve then $\cN_n(X)$ is nonsingular for $n$ large enough (more about this later). This suggests a question, or conjecture: 

{\bf Conjecture 1. } {\it If $X$ is an algebraic variety of any dimension over an algebraically closed field of characteristic zero, then there is an integer $n_0$ such that for $n \ge n_0 $, $\cN _n(X)$ is nonsingular.}

Actually, Yasuda conjectures a more general statement, which implies the one just made, namely:

{\bf Conjecture 2. } {\it  Let char($k$) be zero, $X$ a $d$-dimensional variety over $k$, $\mathcal J _X$ the Jacobian ideal, and $Y$ the closed subscheme of $X$ defined by $\mathcal J _X ^d$. Let $[Z] \in {\cN}_n (X)$ be such that $Z \nsubseteq Y$. Then ${\cN}_n (X) $ is smooth at $[Z]$.}

In fact, it can be proved that for $n$ large enough, every $[Z] \in {\cN}_n (X)$ satisfies $Z \nsubseteq Y$, hence an affirmative answer to Conjecture 2 implies an affirmative answer to Conjecture 1.

Were Conjecture 1 true, then (for $n$ large enough) the $n$-th Nash modification would provide a one-step procedure to resolve the singularities of an algebraic variety (in characteristic zero). Unfortunately, this is not true.    
There are surfaces $X$ such that ${\cN}_n (X)$ is singular for any $n$. We'll discuss such an example, due to Toh-Yama, in Section 7.

There is alternative way to define $\pi _n$, as a special case of the theory of \cite{OZ}. 

Let $\cI$ be the ideal sheaf defining the diagonal 
$\Delta \subset X \times X$. Set $\cP^n _X:=\cO _{X\times X}/\cI ^{n+1}$ and 
 $\cP^n _{X,+}:=\cI _{X\times X}/\cI ^{n+1}          $. These are coherent $\cO_X$-modules, via the first projection. Actually,
  $\cP^n _{X}:=     \cO_X  \oplus     \cP^n _{X,+}         $. 
  In \ref{d3} we considered the affine, or local, version of this construction.  Then it can be shown that $\cN_n(X) $ (as defined above using the Hilbert scheme) can be identified to either 
  $\cN(X,\cP^n _{X})$ or $\cN(X,\cP^n _{X,+})$, in the sense of \cite{OZ}.
  
  This approach (using \cite{OZ}) is significant, because it can be easily adapted to the ``formal'' case of complete local rings, with minor changes, like replacing the tensor product by the completed one. That is, in particular, we may define $\cN_n(Y)$ for a scheme
 $Y={\mathrm {Spec}}(A)$, for a complete local ring $A$, containing $k$ as a field of coefficients.

 In \cite{Ya} it is proved that if $X$ is an algebraic variety and $\widehat X = {\mathrm {Spec}}(\widehat{\cO_{X,x}})$, $x \in X$,  then 
 $$ \cN_n(\widehat X)=\cN_n(X) \times _X \widehat X \, .  $$
 
 Yasuda uses this formal-local theory to investigate properties of $\pi_n$ in the one-dimensional case.
 
 So, now we assume $X$ is a curve. Then in \cite{Ya} it is shown that to prove Conjecture 2 one is reduced to the case 
 $$ (1) \qquad X={\mathrm {Spec} } (R) \,,$$
 where $R$ is the completion of the local ring of a curve at a closed point, and that we may assume $R$ is an integral domain (i.e., what sometimes is called a branch).
 
 So, in this case we have to prove: 
 
 \begin{pro}
 \label{yas}
  If $\cJ$ is the Jacobian ideal of $X$ (or of $R$),  $J \subset X$ is defined by $\cJ$) and $Z$ (corresponding to a point $[Z] \in \cN_n(X)$) is a subscheme of $X$ such that $Z$ is not (scheme-theoretically) contained in $J$, then $[Z]$ is a nonsingular (or normal) point of $\cN_n(X)$.
\end{pro} 
 Note that in this situation (since the normalization $\tilde X$ of $X$ dominates $\cN_n(X)$), 
 $\pi_n:\cN_n(X) \to X$ is a homeomorphism. The only possible singularity of $X$ is the ``origin'' $0$ (corresponding to the maximal ideal of $R$), the only possible singularity of $\cN_n(X)$ is  $0_n$, the only point such that $\pi_n(0_n)=0$. So, in this case $Z$ is a subscheme of $X$ concentrated at $0$, corresponding, say, to an ideal $\cA_n \subset R$; also $[Z]=0_n$. 
 
 The strategy of Yasuda is to use the semigroup (or monoid) 
 $S=S(R)=\{{\mathrm {ord}}(f): f \in R\} = \{s_{-1}, s_0, \ldots, s_j, \ldots\} \subset \mathbb N _0$. 
 Yasuda proves: 

 \medskip
 {\it $0 \in \cN_n(X)$ is nonsingular if and only if $s_n - 1 \in S$. }

\medskip 
 This is a hard result (Theorem 3.3 in \cite{Ya}). Aside from a technical lemma on matrices, the proof uses a criterion for nonsingularity of $0_n$ (\cite{Ya}, Theorem 3.1). To state it, one introduces for $0_n$ an associated embedded first order deformation ${\mathcal Z}_{n,\epsilon} \subset X \times _k  \Spec (k[\epsilon])$, 
 $\epsilon ^2 = 0$. One also has $\mathcal A_n$, the ideal of $R$ defining $Z \subset X$.  Then, the criterion says that the following are equivalent:
 
 (a) $0_n$ is a nonsingular point of $\cN_n(X)$.
 
 (b) ${\mathcal Z}_{n,\epsilon} $ is not the trivial embedded deformation.
 
 (c) Let $k[[y]]$ be the normalization of $R$, $\nu: \tilde X = \Spec (k[[y]]) \to X$ the normalization morphism, 
 $\Gamma _{\nu} \subset \tilde X \hat \times  X =\Spec(R[[y]])    $ the graph of $\nu$, $P$ be the prime ideal of $R[[y]]$ defining $\Gamma _{\nu}$ and 
 $P^{(n+1)}$ the $(n+1)$-th symbolic power of $P$. Then (c) states: there is a series $f = a_0 + a_1 y + \cdots $ in  $P^{(n+1)} \subset R[[y]]$ such that 
 $a_1 \notin {\mathcal A }_n$.
 
 The proof of \ref{yas} also uses the conductor $\mathcal C$ of $R$ in $k[[y]]$, or the subscheme $C $ of $X$ it defines. The ideal $\mathcal C$ of $R$  
  is generated by monomials $y^j$, $j \ge t_l$ for a suitable integer $t_l$.  We have $\cJ \subseteq \cC$. Also, in the proof of 3.3.in \cite{Ya} it is proved that if 
  $s_n \le t_l+1$ then $\cC \subset \cA_n$, i.e., in this case  we have 
   inclusions of  subschemes of $X$:
   $$Z \subseteq C \subseteq J \,.$$
   With these preliminary results, the proof of Conjecture 2 is easy. Indeed assume $Z \subseteq J$. Consider 
   $S=S(R) = \{ s_{-1}, s_0, \ldots, s_n, \ldots  \}$. We cannot have $s_{n} - 1< t_l$, because then 
    $ \subseteq C \subseteq J$, contrary to the assumption $Z \nsubseteq J$ of our conjecture. Thus, 
    $s_n - 1 \ge t_l$. But then $s_n - 1 \in S$, because any $j \ge t_l$ is in $S$. By 3.3 of \cite{Ya}, $0_n$ is nonsingular.

A similar result is false in positive characteristic. In fact, Yasuda proves that if $X=\Spec(R)$ is a branch as in the previous discussion, but now the characteristic of the base field is $p > 0$, then for any natural number $e$ large enough $\cN_{p^e-1}(X) \cong X$. 

\section{The normal case}
\label{S:nor}
We have seen that the basic smoothness result (Theorem \ref{basic}) fails in positive characteristic. The simplest counterexample involves a singular curve whose Nash modification is an isomorphism. This variety, being one-dimensional, is also non normal. The non normality is essential in this example, because we have the following result.

\begin{thm} Let $X$ be a normal algebraic variety over a field $k$, algebraically closed of characteristic $p > 0$. Then, the Nash modification $\pi : \cN(X) \to X$ is an isomorphism if and only 
if $X$ is singular.
\end{thm} 
This is due to Duarte and Betancourt (\cite{DNB}). We sketch their proof, recalling the basic necessary facts. One implication is clear, the hard one is to see that if $\pi$ is an isomorphism then $X$ is non singular. So, let $x$ be a (closed) point of $X$, $R=\cO_{X,x}$ its local ring. We shall see that if $\pi$ is an isomorphism then $R$ is a regular local ring. 

They use a result of E. Kunz (\cite{K}). Before we state it, we review some terminology. Let $R$ be a domain which is  an  algebra over a field  $k$  of positive characteristic $p$, with field of fractions $F$, and $L$ an algebraic closure of   $F$. Then, 
$R^{1/p}=\{a \in L: a^p \in R \}$. We have $R \subseteq R^{1/p}$. 

Kunz showed that if $R^{1/p}$ is a finite free $R$-module (i.e., $R^{1/p}$ is isomorphic to $R^m$, for some integer $m$) then $R$ is regular. 

Duarte and Betancourt find a free $R$-module with free basis $\{ e_{\alpha}: \alpha \in \cA\}$, where $\cA$ is a finite set, and a homomorphism 
$\psi : \sum_{\alpha \in \cA} Re_{\alpha} \to R^{1/p}$, and they show that $\psi$ is bijective. In the proof of bijectivity, they use the following fact: if $\psi$ induces a bijection
$$\psi_Q : \sum_{\alpha \in \cA} R_Q e_{\alpha} \to R_Q^{1/p}$$
whenever $Q=0$ (the zero ideal) or dim$(R_Q) =1$, then $\psi$ is bijective.

It is here that the assumption on normality plays a role: if $R$ is normal then each localization $R_Q$ of dimension one is a discrete valuation ring, this allows us to use special arguments.

Next we provide more details.

The set $\cA$ is $\{ \alpha=(\alpha_1, \ldots, \alpha_d )\in {\mathbb N}^d: 0 \le \alpha_i <p  \}$, it has $p^d$ elements. We claim:

(a) from the fact that $\pi: \cN(X) \to X$ is an isomorphism, we can find a system of parameters $x_1, \ldots, x_d$ of $R$, and $k$-derivations $\delta _1, \ldots, \delta_d$, $\delta _i:R \to R$, $i=1, \ldots, d$, such that the $d$ by $ d$ matrix $(\delta_i(x_j))$, $i,j = i, \ldots d$, is the identity matrix. This will be proved later.

(b) For $\alpha =(\alpha_1, \ldots, \alpha _d) \in \cA$, let 
$$ \delta^{(\alpha)}=   \frac{1}{ {\alpha_1}! \cdots  {\alpha_d}! } \delta^{(\alpha_1)} \ldots \delta^{(\alpha_d)}   $$ 
Note that 
 $\frac{1}{ {\alpha_1}! \cdots  {\alpha_d}! } \in k$, because $0 \le \alpha _i < p$, for all $i$. We claim that there are coefficients ${\tilde c}_{\gamma,\alpha} \in k$, for each $\alpha, \gamma$ in $\cA$, such that for each $\gamma$ the function 
 $\Phi_{\gamma}= \sum_{\alpha \in \cA}  {\tilde c}_{\gamma,\alpha} \delta^{\alpha}$ satisfies: 
 $\Phi_{\gamma}(x^{\beta})=0$ if $\gamma \not= \beta$ and 
$ \Phi_{\gamma}(x^{\gamma})=1$. The maps $\Phi_{\gamma}$ are additive and $k$-linear, but not necessarily $R$-linear. But using the fact that each $\Phi_{\gamma}$ is a composition of derivations and the fact that $p.1=0$ in $k$, $\Phi_{\gamma}$
 is $R^p$-linear (where $R^p=\{b^p:b \in R\} \subseteq R$. Hence, if $\phi_{\gamma}:R^{1/p} \to R^{1/p}$ is defined by 
 $$\phi_{\gamma}(f^{1/p})=[\Phi_{\gamma}(f)]^{1/p} ~,$$
 we have $\phi_{\gamma} \in \mathrm {Hom} _{R}(R^{1/p}, R^{1/p})$.  
 
 Now, if for $\beta = (\beta_1, \ldots, \beta_d) \in {\mathbb N}^d$, $x=(x_1, \ldots, x_d)$ we write 
 $$x^{\beta/p}=((x_1^{1/p})^{\beta_1}, \ldots, (x_d^{1/p})^{\beta_d}) \, ,$$ 
 we have
 $\phi_{\gamma}(x^{\gamma /p})=1$ and  $\phi_{\gamma}(x^{\beta /p})=0$ if $\gamma \not= \beta$.
 
(c) Now we may explain the map
 $$(\star) \qquad \psi: \bigoplus _{\alpha \in \cA}R{e_{\alpha}} \to R^{1/p}$$
mentioned before. It is enough to describe how it operates on elements of the basis $\{e_{\alpha}\}$ :
$$\psi (e_{\alpha})=x^{\alpha/p} \in R^{1/p} ~.$$
As remarked, to show the bijectivity of $\psi$ it suffices to show the bijectivity of the localization 
$\psi_{Q}$ at each prime ideal $Q$ of $R$ such that $\mathrm{dim} R_Q \le 1$.

First of all, note that from work of E. Kunz and K. Smith, since $R$ is excellent (being the local ring of an algebraic variety) and $k$ is perfect (being algebraically closed), $R^{1/p}$ is a finite $R$-module. The same is true for any localization at a prime. Moreover, the rank of $R$ (i.e., the dimension of the $F$-vector space $R^{1/p} \otimes _R F$, $F=\mathrm{Frac} (R)$), is equal to $p^d$. 

Now, if $\mathrm{dim} (R_Q)=1$, then $R_Q^{1/p}$ is free. This is because, by the normality of $R$, $R_Q$ is regular, even a discrete valuation ring (hence a principal ideal domain). Then the assertion follows from the cited theorem of Kunz or, alternatively, because $R_Q^{1/p}$ is finitely generated and torsion free (since $R_Q^{1/p} \subset F$) and $R_Q$ is a principal ideal domain. The rank of $R_Q^{1/p}$  is again $p^d$.

Hence the inclusion $j:R_Q \to R_Q^{1/p}$ admits a section  
$\sigma:  {R_Q}^{1'p} \to R_Q$, so that $\sigma j =id_{R_Q}$. Let 
$$\rho _Q: {R_Q}^{1/p} \to \bigoplus _{\alpha \in \cA}{R_Q}{e_{\alpha}} $$
be defined by 
$\rho_Q({f^{1/p}}/{s})=\sum_{\alpha \in \cA}\sigma({\phi_{\alpha}(f^{1/p})}/{s})$, for $f \in R$, $s \in R \setminus Q$.

Then one verifies that $\rho_Q(x^{\alpha/p})=1.e_{\alpha}$, whence $\rho _Q$ is onto. Since both the domain and codomain of $\rho$ are free $R_Q$-modules of the same rank ($=p^d$) it follows that $\rho _Q$ is also injective, i.e., an isomorphism. Also, $\rho _Q \Psi_Q$ is the identity of 
$\bigoplus _{\alpha \in \cA}R{e_{\alpha}}$. It follows that $\Psi_Q$ is the inverse of $\rho_Q$, hence also an isomorphism, as needed. 

The case where $\mathrm{dim}(R_Q)=0$ is similar, but simpler, because now $R_Q$ is the field $F$, the fraction field of $R$. 

(d) Now we justify the claim made in (a) regarding the parameters $x_i$ and the derivations $\delta _j$. We have to recall some ``basic'' facts about differential operators and the module of principal parts. What follows can be developed with greater generality, but to simplify we work with an algebraically closed field $k$. 

(i) If $R$ is a (noetherian) $k$-algebra, let  $I$ be the kernel of the multiplication homomorphism 
$R \otimes _k R \to R$ and consider $R/I^{n+1}$, which is an $R$-module, via the rule
 $r(a \otimes b)= ra \otimes b$). This is the {\it module of principal parts of order n} of $R$, denoted by 
 $\cP^n_{R}$.
 
 Then one can prove rather easily, using the universal property of Grothendieck's Grassmanian, cf. Section 3 that (using the notation introduced at the beginning of this proof), if 
 $\pi : \cN(X) \to X $ is an isomorphism then 
 $$\cP^n_R=R^{\binom{n+d}{d}} \oplus T ~ ,$$
 where $T$ is a torsion $R$-module.
 
 (ii) {\it Differential operators of order $\le n$, over $k$}. These are the $R$-modules $D^n_{R|k}$, 
 $n=0, 1, 2, \ldots $ defined inductively as follows:
 
 $D^0_{R|k}= {\mathrm {Hom}}_R(R,R)$ 
 
 If $D^n_{R|k}$ was defined,    $D^{n+1}_{R|k}=
 \{ \delta \in {\mathrm {Hom}}_k(R,R): \delta r - r \delta \in          D^n_{R|k} \}$.
 
 Then $D^{n}_{R|k} \subseteq D^{n+1}_{R|k}$, for all $n$. The union of all the $D^n_{R|k}$ is naturally a ring (by composition), denoted by $D_{R|k}$
 
(iii) {\it Derivations} $Der_{R|k}$. We let $Der_{R|k}$ denote the $R$-module of $k$-derivations of the $k$-algebra $R$. 

It can be proved that $D^{1}_{R|k}  = R \oplus Der_k(R)$.

(iv) {\it Differential powers}. Let $R$ be a local $k$-algebra, with maximal ideal $M$, and $A/M=k$. The {\it n-th differential power of M} is the $R$-ideal  
$M^{\langle n \rangle} = \{ a \in R : \delta (a) \in M, ~ \forall \, \delta \in D^{n-1}_{R|k}  \}$

It can be proved that 
$$   \mathrm{dim} (R/ M^{\langle n+1 \rangle})=f.r.(\cP^n_{R|k})  ~,       $$
where dim denotes dimension as a vector space and $f.r.$ denotes free rank, that is the largest rank of a free $R$-module which is a quotient of $\cP^n_{R|k}$.

(v) {\it A perfect pairing}. With the notation above, there is a perfect pairing of $k$-vector spaces:
$$(1) \quad   D^{n-1}_{R|k} / \cI _n \times R/M^{\langle n \rangle} \to k=R/M  ~,           $$
where $\cI _n = \{ \delta \in   D^{n-1}_{R|k}: \delta(R) \subseteq M          \}$, given by the rule 
$({\bar \delta} , {\bar r}) \mapsto \bar{\delta (r)} $.

Of special interest to us is the case $n=2$. In this case, using the equality  
$D^{1}_{R|k}  = R \oplus Der_k(R)$ of (iii), (1)  becomes 
$$ (2) \quad  Der_k(R)/ \cI'_n \times  M/    M^{\langle 2 \rangle}     \to k  ~ ,$$
where  $\cI' _n = \{ \alpha \in      Der_k(R)  : \alpha(R) \subseteq M          \}$.
 
 Now, by (iv) we know that 
$\mathrm {dim}(R/M^{\langle 2 \rangle}) 	=    f.r. (\cP^1_{R|k})   $. But by (1) in (i), 
$   \cP^1_{R|k} = R^{d+1}\oplus T           $, with $T$ torsion. So, $ f.r. (\cP^1_{R|k})=d+1$. It follows that 
      $ \mathrm {dim}(M/M^{\langle 2 \rangle}) =        \mathrm {dim}(R/M^{\langle 2 \rangle}) -1 = d$.
      
 Now we are in position to justify the claim in (a) in this section. Choose a basis $\bar{x_1}, \ldots, \bar{x_d}$ of 
 $M/M^{\langle 2 \rangle}$, and (by the pairing (2) of (v)) the dual basis 
 $\bar{\partial _1}, \ldots, \bar{\partial _d}$. Lift them to $M$ and $Der_k(R)$, getting 
 $x_1, \ldots, x_d$ and  ${\partial _1}, \ldots, {\partial _d}$ respectively. The, 
 $\partial_i (x_j)$ is a unit of $R$ if and only if $i=j$.
 
 The derivations $\partial_i$ do not necessarily have all the properties of the claimed derivations $\delta_i$ of (a). But we may modify them as follows. Let $A=(a_{ij})$ be the $d \times d$ matrix whose $ij$ coefficient id $\partial_i(x_j)$. This matrix is invertible, since when we reduce the coefficients mod $M$ we get the identity. Let $C=(c_{ij})$ be its inverse in $GL_d(R)$, and 
 $\delta_k = \sum_{i=1}^{d} c_{ki} \partial_i .$
 Then,
 $$\delta_t(x_j)=\sum_{i=1}^{d} c_{ki} \partial_i =        \sum_{i=1}^{d} c_{ti} a_{ij} ~.$$
 Hence $\delta_j(x_j) =1$ and $\delta_t(x_j) = 0$ for $t \not= j$. Thus 
 $\delta_j, \, j=1, \ldots, d$ are the claimed derivations.  
 
 \section{Toric varieties}
 \label{S:toro}
 
 In case our variety $X$ is toric, Theorem \ref{basic} can be strengthened to a statement involving the $n$-th Nash modification. Namely, we have the following result. As before, in this section the base field $k$ is algebraically closed. 
\begin{thm}
\label{tor}
Let $X$ be a normal toric variety over $k$, of arbitrary characteristic, $n$ any positive integer, $\pi_n : \cN_n(X) \to X$ the $n$-th Nash modification. Then $\pi_n$ is an isomorphism if and only if $X$ is an nonsingular. 
\end{thm}
This result follows from: 
\begin{thm}
\label{torno}
Let $X$ be a normal toric variety over $k$, $n$ any positive integer, $\bar{\pi}_n : \cN^-_n(X) \to X$ the normalized $n$-th Nash modification. Then $\bar{\pi}_n$ is an isomorphism if and only if $X$ is an isomorphism. 
\end{thm}

The hypothesis on normality in these statements is really not necessary, since a toric variety (according to the definition we use) is automatically normal.

To see that Theorem \ref{torno} implies Theorem \ref{tor} , note that if $\bar{\pi}_n $ is an isomorphism then, since $X$ is normal,  $\cN(X)$ is also normal normal. Hence 
the normalization map 
$\eta:\cN^-(X) \to \cN(X)$ is an isomorphism and the implication follows.

As usual, we may assume that $X$ is affine. These theorems are proved in \cite{D} in characteristic zero and in \cite{DB} in positive characteristic. We'll sketch a proof of Theorem \ref{torno}.

We begin by recalling some basic notions of toric geometry. Usually we omit proofs, since there are many good references on toric varieties, e.g.,  \cite{C} \cite{F}, \cite{M},  \cite{O}, as well as many notes available on the web.

\medskip

(a) {\it Toric varieties.} Let $N$ be a lattice, that is a finitely generated free abelian group, $N \cong {\mathbb Z}^d$. We let $N_{\mathbb R}=N \otimes _{\mathbb Z}{\mathbb R} \cong {\mathbb R}^d$, 
$M=N^*= \mathrm{Hom} _{\mathbb Z}(N, \mathbb Z)$, 
$M_{\mathbb R}=M \otimes _{\mathbb Z}{\mathbb R}$, isomorphic to the dual of $N_{\mathbb R}$

A strongly convex, rational, polyhedral  cone (SCRPC) is a subset $\sigma \subseteq N_{\mathbb R}$ of the form 
$\sigma=\bR_{+}v_1 + \cdots + \bR_{+}v_s$, where $\bR_{+}$ denotes the set of nonnegative reals, and $v_i \in N , ~ {\forall i}$; we require that $\sigma $ does not contain any line. In general, a SCRPC will be simply called ``a cone''. 

If $\sigma$ is a cone, $\sigma ^{\vee}= \{ u \in M_{\bR}: u(v) \ge 0, ~ \forall v \in \sigma  \}$, which is a cone in $M_{\bR}$. 

The {\it semigroup} of a cone $\sigma$ is 
$S_{\sigma}= \sigma^{\vee} \cap M  =\{ u \in M: u_{|\sigma} \ge 0  \}$. This is a finitely generated semigroup or, actually, a monoid (Gordan's Lemma).

Let $k[S_\sigma]$ be the $k$-semigroup algebra corresponding to $\sigma$. That is we consider a symbol $\chi^s$ for each $s \in S_{\sigma}$, then elements of $k[S_\sigma]$ are finite linear combinations $\sum_i r_i \chi^{s_i}$,  $r_i \in k, ~ s_i \in S_{\sigma}$, $\forall i$; the product is according to the rule: $\chi^s . \chi^t =\chi^{s+t}$.  

Assume $k$ is an algebraically closed field. The scheme $X=U_{\sigma}=\Spec(k[S_{\sigma}])$ is the {\it toric affine variety } defined by the cone $\sigma$. As usual, in general we shall work with $X(k)$, the closed points of $X$ only. We have a natural bijection 
$$X(k)= {\mathrm{Hom}} _k(k[S_{\sigma}], k)=      {\mathrm{Hom}} _{sg}(S_{\sigma},k) ~,$$
where, in ``${\mathrm{Hom}} _{sg}(S_{\sigma},k)$'',  $k$ is regarded as a semigroup under addition.

The variety $S_{\sigma}$ is irreducible and normal (see \cite{F}, page 29), but not necessarily regular.  A necessary and sufficient condition for $S_{\sigma}$ to be regular is that  the cone $\sigma$ be equal to $\bR_{+} e_1 + \cdots + \bR_{+} e_s$, where $e_1, \ldots, e_s$ (which are elements of $N$) are part of  
 a basis of $N$ (as a free abelian group), see \cite{F}, page 29.) 

A face of $\sigma $ is a subset $\tau \subset \sigma$ of the form $\{x \in \sigma: u(x)=0\}$, for a fixed $u \in M_{\bR}$.

Faces have many nice properties. We list some useful ones:  (i) A face of a cone is a cone (as usual, ``cone'' means SCRPC).  (ii) A face of a face of a cone $\sigma$ is a face of $\sigma$. (iii) The intersection of faces of $\sigma$ is a face of $\sigma$. (iv) If $\tau$ is a face of $\sigma$, defined by $u \in S_{\sigma}$, then $S_{\tau}=S_{\sigma} + {\mathbb N}(-u)$ (hence $S_{\sigma} \subset S_{\tau}$), and the induced morphism $S_{\tau} \to S_{\sigma}$ is an open immersion.

 The origin $\{0\}$ is always a face of a cone $\sigma$. In view of (iv), there is an immersion of 
$\mathbb T_d = S_{0}$ (the $d$-torus, isomorphic to $G_d$ or $(k^{\star})^d$) into $S_{\sigma}$, as a dense open set. 

\medskip

{\it A remark on normality}. Affine toric varieties, in the sense described above, are automatically normal. There is a more general concept of affine toric variety, that yields examples of non normal ones. Namely, if $k$ is an algebraically closed field (of any characteristic) define {\it affine toric variety} (in $k^n$, or $A^n$) as the closed image $X$ of a monomial map from the torus $T_n ={k^{\star}}^n $ to $k^n$ of the form:
$$ \mathbf{x} \mapsto  ({\mathbf{x}}^{a_1}  , \ldots  {\mathbf{x}}^{a_n}    )$$
where $a_i = (a_{i1}, \ldots, a_{in}) \in {\bZ}^n$ are given, $  \mathbf{x} = (x_1, \ldots, x_n  )$ and 
$ \mathbf{x}^{a_i} = ({x_1}^{a_{i1}}, \ldots,     {x_n}^{a_{in}} )  $.

It can be proved that such an $X$ contains a copy of of the $d$-torus $T_d$ (for a suitable $d$), as  dense set, $d$ is the dimension of $X$. Affine varieties as we introduced them are a special case of this notion, namely the normal ones. See \cite{MM} or, for a different approach,\cite{GT}.

In this article we shall be concerned exclusively with normal toric varieties.

{\it Fans and general toric varieties}. A fan $\Sigma$ in $N_{\bR}$ is a finite nonempty collection of cones in $N_{\bR}$ such that: (a) all faces of a cone in $\Sigma$ are in $\Sigma$, (b) if $\sigma, \, \sigma'$ are in $\Sigma$, then $\sigma \cap \sigma'$ is a face of both $\sigma$ and $\sigma'$.

The different varieties $U_{\sigma}$, for $\sigma$ in a fan $ \Sigma$, glue together to produce a variety 
 $X=X_{\Sigma}$, which is integral and normal, called a toric variety. If dim($X =d$), the the $d$-torus $\mathbb T _d$ can be identified to an open set in $X$, in such a way that the natural action of $\mathbb T _d$ on itself (by multiplication) extends to an action of $\mathbb T _d$ on $X$.
 
 Conversely, it can be proved that an integral, normal variety $X$  of dimension $d$, containing an open set isomorphic to the $d$-torus, such that the action of $\mathbb T _d$ on itself  extends to an action of $\mathbb T _d$ on $X$ is isomorphic to $X_{\Sigma}$,  in a $\mathbb T_d$-equivariant way, for a suitable fan $\Sigma$.
 
 Suppose $\Sigma, \, \Sigma'$ are fans in $N_{\bR}$ such that every cone in  $\Sigma$ is contained in a cone in $\Sigma'$ and every cone in $\Sigma'$ is union of cones in $\Sigma$. Then we say $\Sigma$ is a refinement of $\Sigma'$. If $\Sigma$ is a refinement of $\Sigma'$ there is an induced morphism 
 $f:X_{\Sigma} \to X_{\Sigma'}$, which is proper, $\mathbb T _d$-equivariant and birational.
 
 \medskip
{\it Groebner fans.} In the proof of Theorem \ref{torno}, an important role is played by a Groebner fan.
 
 This is a fan associated to an ideal in the ring $k[x_1, \ldots, x_d]$ or, more generally, to a ring 
 $R \subseteq k[x_1, \ldots, x_d]$ of the following form. Consider $d$-tuples 
 $a_i=(a_{i1}, \ldots, a_{id}) \in {\mathbb N}^d$, $i=1, \ldots, s$. Then 
 $R=k[x^{a_1}, \ldots, x^{a_s}]$, where 
 $x^{a_i}={x_1}^{a_{i1} }\ldots  {x_d}^{a_{id}}$. 
 
 If $<$ is a well-ordering in ${\mathbb N}^d$, it induces a well-ordering in the monomials $x^{\beta}$ in $R$, namely $x^{\beta} < x^{\beta'}$ if $\beta < \beta'$. If for  $\alpha, \beta, \gamma$ in  ${\mathbb N}^d$, from     $\alpha <\beta$ it follows $\alpha + \gamma <\beta + \gamma$, we say that we have {\it a monomial ordering} in $R$.
 
 If $w \in \bR^d$ and $f = \sum c_ux^u$, $c_u \in k$ $\forall u$, is an element  of $R$, we may associate to $f$ a polynomial 
 $in_w(f) \in k[x_1, \ldots, x_d]$, the {\it initial form of  $f$ relative to $w$}, namely
  the sum of terms $c_u x^u$ such that $w.u$  (= $\sum w_i u_i$) is a maximum. Sometimes we simply write $wu$, rather than $w.u$. 
  
  If $I$ is an ideal of $R$, $in_w(I)$ is the ideal generated by $\{  in_w(f): f \in I \}$.
  
 Next assume $\sigma$ is a cone with $\sigma^{\vee}={\bR}_{+} (a_1, \ldots, a_s) \subset {\bR}_{+}^d$. 
 We also assume $N$, the lattice of $\sigma$, is ${\mathbb Z}^d$, so that the duality  $N \times M \to {\mathbb Z}$ is the usual scalar product. With notation as above, for $w \in \sigma$ let 
 $C[w]=\{ z \in \sigma: in_w(I) =in_z(I)\}$. Then it can be proved that 
 ${\overline C}[w]$,   
  the closure  (in $\bR^d$)  of $C[w]$ is a cone  (i.e., a SCRPC), $C[w]$ is the relative interior of ${\overline C}[w]$ and, as $w$ varies in $\sigma$, the cones ${\overline C}[w]$ form a fan, the so-called Groebner fan of $I$. We'll denote it by $GF(I)$.
  
 In the verification of the fact that  ${\overline C}[w]$ is a cone, whose interior is $C[w]$, one uses an auxiliary monomial order in $R$. Namely, pick any monomial ordering $\succ$ in $R$ and let define 
 $x^{u}>x^v$ if either $uw > vw$ or , if $uw=vw$, require $x^u  \succ x^v$.
 
 In the proof of Theorem \ref{torno}, we may assume $X$ is the affine toric variety $U_{\sigma}$, where $\sigma$ is such that $S_{\sigma}=k[\sigma^{\vee}\cap M]=k[x^{a_1} , \ldots, x^{a_s} ]=R$, with $a_i \in {\mathbb N}^d$, $\forall i$.
 
 It is easy to see that $\cN_n^{-}(X)$, the normalized $n$-th Nash modification of the affine toric variety $X$,  is again a toric variety. Indeed it contains a copy of the $d$-torus 
 $\mathbb T _d$, acting on it as necessary. Moreover, the morphism ${\bar \pi}_n$ is $\mathbb T _d$-equivariant. It is a general fact that then $\cN_n^{-}(X)$ is isomorphic to $X_{\Sigma}$, where $\Sigma$ is a refinement of the trivial fan consisting of all the faces of $\sigma$. A key element of the proof of \ref{torno} is to verify that we may assume that $\Sigma$ is the Groebner fan of a certain ideal of $R$, namely 
 $I_n=\langle x^{a_1} -1, \cdots,  x^{a_s} -1  \rangle$. 
 
 To check that $\Sigma = GF(I_n)$ is not easy. Among other things, one uses the notions of distinguished point of a cone and of limits. If $\sigma$ is a cone, its {\it distinguished point} is a closed, or $k$-rational, point  $x_{\sigma}$ obtained as follows. Consider the following homomorphism $\phi$ of semigroups
 $S_{\sigma} \to k$, sending a $u \in S_{\sigma} \subseteq M$ to $1 \in k$ if $u_{| \sigma} =0$, or to $0$ if $u(w) >0$ for some $w \in \sigma$. To $\phi$ it corresponds a unique morphism 
 $\Spec (k) \to U_{\sigma}=X$, this element of $X(k)$ is $x_{\sigma}$. This distinguished point may be introcuced in a more ``dynamic'' way using one-parameters subgroups of ${\mathbb T}_d$, i.e., 
 homomorphisms of algebraic groups $G_m \to {\mathbb T}_d$. Identifying $G_m=k^{\star}$ and 
 ${\mathbb T}_d = (k^{\star})^d$, such a $\lambda$ is of the form 
 $\lambda (z) = (z^{c_1},  \ldots, z^{c_d})$, each $c_i \in \mathbb Z$ (perhaps $c_i < 0$ for some indices). Since ${\mathbb T}_d \subset U_{\sigma}$, $\lambda$ may be regarded as a map into $U_{\sigma}. $Then it can be proved that if $v $ is in the interior of $\sigma$, 
 $\lim _{z \to 0} \lambda (z) = x_{\sigma}$,
  if $v$ is in the interior of a face $\gamma$ of $\sigma$,  
  $\lim _{z \to 0} \lambda (z) = x_{\gamma}$, and if $v \notin \sigma$, then the limit does not exist. If $k = \mathbb C$, these are limits in the  sense of Calculus, for more general fields they must be described more algebraically.  See \cite{F} or \cite{M}.
  
  Once we know that $\cN _n ^{-}(X)$ is the toric variety associated to the Groebner fan $GF(I_n)$, to show that if ${\bar \pi} _n$ is not an isomorphism then $U_{\sigma}$ is not regular, it suffices to see that the Groebner fan $GF (I_n)$ is not ``trivial'', i.e., it is not the fan consisting of all the 
  faces of $\sigma$. 
  To verify this it is sufficient to find points $w, \, w'$ in the interior of 
  $\sigma$  so that 
  $in_w (I_n) \not= in_{w'} (I_n)$. Actually, it is enough to verify a statement that not involve the full ideal $I_n$. Namely, consider in $R=k[x^{a_1}, \ldots, x^{a_s}]$ the order $>_{w}$ mentioned before (when we introduced the Groebner fan) and the reduced Groebner basis $\cB$ of $R$ relative to this order. Then to show the inequality $in_w (I_n) \not= in_{w'} (I_n)$ it suffices to find $g \in \cB$ such that 
  $in_w (g) \not= in_{w'} (g)$. 
  
  These steps are accomplished if $char (k) =0$ in \cite{D}. If   $char (k) > 0$, in \cite{DB} it is proved that the arguments of \cite{D}, conveniently modified, are still valid in this case.
  
 In dimension 2, more precise results on resolution of singularities are available. Namely, if $X$ is a toric surface over a field $k$ of chacteristic zero, G. Gonzalez Sprinberg showed that a finite sequence of normalized Nash transformations desingularizes $X$ (see \cite{GS}). Recently, Daniel Duarte, Jack Jeffries and Luis Nu{\~n}ez-Betancourt proved that the result is also valid in positive characteristic (see \cite{DJN}).  
  
   \section{Toh-Yama's example}
   \label{S:Toh}
   
   In \cite{Ya}, Yasuda proved that if $X$ is an algebraic curve over an algebraically closed field of zero characteristic, then $\cN_n(X)$ is nonsingular for a suitable $n \in {\mathbb N}$. It is also conjectured that a similar statement might be true for $X$ a variety of arbitrary dimension. We reviewed this work in Section 4. 
   
   But in the paper \cite{TY}, Toh-Yama gives an example showing that this is not true. Indeed, if $X$ is the surface (in 
   ${\mathbb A}_{\mathbb C}^3 = \Spec ({\bC}[x,y,z] $) defined by $z^4-xy=0$, then 
   $\cN_n^{-}(X)$ is singular for all $n$. From this it easily follows that $\cN_n(X)$ is also singular for all 
   $n$. The surface $X$ has a single (rational) singularity at the origin, called the $A_3$-singularity. In \cite{TY} one works over $\bC$ and methods of toric geometry are used. We sketch the main lines of the proof.
   
The surface $X$ is toric: $X=U_{\sigma}$, the affine toric variety corresponding to the cone  
$\sigma \subset {\bR}^2$ generated by $(0,1)$ and $(4,-3)$. The semigroup $S_{\sigma}$ is generated by $(1,0), (3,4)$  and $(1,1)$, and we may write ${\bC}[S_{\sigma}]={\bC}[u, u^3v^4,uv]$. As explained in Section \ref{S:toro}, $\cN_n^{-}(X)$ is the toric variety associated to the Groebner fan $GF(J_n)$, where 
$J_n = \langle u-1, u^3v^4-1,uv -1\rangle$. But Toh-Yama introduces another, more combinatorial,  technique to produce cones of $GF(J_n)$. Let us explain it, in a more general setting. 

Working over an algebraically closed base field $k$, possibly of positive characteristic, consider a cone $\sigma \subset \bR^d$ whose semigroup ring 
$R=k[S_{\sigma}]$ 
is of the form $R = k[x^{a_1}, \ldots, x^{a_n}]$ (notation as in the previous section, i.e., $x=(x_1, \ldots, x_d)$, $a_i \in {\bN}^d, \forall i$). Assume a monomial ordering $\prec$ is defined on $S$.  A 
{\it marked reduced Groebner basis} (MRGB)  of $I$ (an ideal of $R$) s a finite collection of ordered pairs 
$\cG =\{ (g_1, \alpha _1)  , \ldots, (g_t, \alpha _t) \}$, where $g_1, \ldots, g_t$ is a reduced Groebner basis of $I$, and $\alpha_i$ is the exponent of the leading monomial of $g_i$, for all $i$. We let 
$lm (\cG) = \{ \alpha_1, \ldots, \alpha_t \}$. 

To such MRGB $\cG$ we associate a set $C_{\cG} \subset \bR^d$, namely 
$$C_{\cG} = \{ w \in \sigma : (\alpha_i - \beta). w \ge 0, \forall \beta \in Sup(g_i), i=1, \ldots, t \}$$
where $Sup(g_i)$ is the set of exponents $\beta$ such that $x^{\beta}$ appears in $g_i$ with nonzero coefficient, and the dot denotes the usual scalar product. This is a (strongly convex, rational polyhedral) cone $\tau \subseteq \sigma$, moreover $\tau \in GF(I)$.

This will be applied to our situation, i.e., $U_{\sigma} = X$, the $A_3$-singularity. We introduce in 
$S={\bC}[u, u^3v^4,uv] $ the monomial order $\prec$ associated to the two by two matrix $\mathcal M$ whose rows are $(2 ~-1)$ and $(1 ~1)$. That is, 
$$  u^{\alpha _1}   v^{\alpha _2} \prec   u^{\beta _1}   v^{\beta _2}   $$
if  either $2 \alpha _1 - \alpha_2 < 2 \beta_1 -\beta_2$ or  $2 \alpha _1 - \alpha_2 =2 \beta_1 -\beta_2$ and 
$\alpha _1 + \alpha _2 < \beta_1 + \beta _2$.

We consider the a marked reduced Groebner basis $\cG_n$ of $J_n$ with this order, let 
$\cG_n =\{ (g_1, \gamma_1), \ldots, (g_s, \gamma _s)$, so 
$lm (\cG_n) =\{ \gamma_1, \ldots, \gamma _s \}$. An important item in Toh-Yama's example is an alternative description of $lm(\cG_n)$. Indeed, in a more formal way a set $\cP_n$ of pairs of numbers is introduced. Although this is not obvious, in a rather complicated way it is checked that $\cP_n=lm(\cG_n)$.

Using $\cP_n$, it is possible to give an explicit description of the (two) rays bounding the cone 
$\tau := C_{\cG_n} \subset \bR^2$. The answer is: one ray is $L_1=\bR_+ (2,-1)$. The other is 
$L_2=\bR_+ l_n$, where
\begin{itemize}
\item $l_n = (2n-2,-n+2$), if $n$ is odd,
\item $l_n = (2n,-n+1$), if $n$ is even.
\end{itemize}

Now it is easy to check that in either case $U_{\tau}$ is singular, that is $\tau$ is not generated by a subset of a basis of a free $\bZ$-module.

In \cite{DB} the authors show that To-Yahma's example is still valid if the base field is of positive characteristic. They use the same method as Toh-Yama's, including the introduction of the sets $\cP_n$ to obtain useful cones in the Groebner fan. The whole proof of \cite{TY} does not directly work in positive characteristic, but they manage to reduce their problem to a situation where the results of that paper apply.

\providecommand{\bysame}{\leavevmode\hbox to3em{\hrulefill}\thinspace}

\end{document}